\documentclass{cmam}
\usepackage{amsmath,amssymb,amsfonts,xcolor}
\usepackage[ansinew]{inputenc}
\usepackage{epsfig,cite}
\usepackage {longtable}

\newcommand{{\un}}{u^{(0)}_{n}(x)}

\newcommand{{\bea}}{\begin{array}}
\newcommand{{\ea}}{\end{array}}
\newcommand{\be}{\begin{equation}}
\newcommand{\ee}{ \end{equation}}

\newcommand{{\il}}{\int}
\newcommand{{\la}}{\lambda_n^{(0)}}

\begin{document}
\runauthor{I.~Gavrilyuk, V.~Makarov, N.~Romaniuk} \runtitle{Eigenpairs of Fractional Jacobi-Type Equation}
\begin{topmatter}

\bibliographystyle{cmam}

\title {Super-Exponentially Convergent Parallel
Algorithm for a Fractional Eigenvalue
Problem of Jacobi-Type}

\author{Ivan Gavrilyuk}
\address{University of Cooperative Education Gera-Eisenach, 
Am Wartenberg 2, 99817 Eisenach, Germany}
\email{iwan.gawriljuk@dhge.de;  http://orcid.org/0000-0003-3115-9690}
\author{Volodymyr Makarov}
\address{Institute of Mathematics of NAS of Ukraine, 3 Tereshchenkivs’ka Str., 01004 Kyiv-4, Ukraine}
\email{makarov@imath.kiev.ua; http://orcid.org/0000-0002-4883-6574}
\author{Nataliia Romaniuk}
\address{Institute of Mathematics of NAS of Ukraine, 3 Tereshchenkivs’ka Str., 01004 Kyiv-4, Ukraine}
\email{romaniuknm@gmail.com; http://orcid.org/0000-0002-3497-7077}

\begin{abstract}
\noindent A new algorithm for eigenvalue problems for the fractional Jacobi type ODE is proposed. The
algorithm  is based on piecewise approximation of the coefficients
of the differential equation with subsequent recursive procedure adapted
from some homotopy considerations. As a result, the  eigenvalue problem (which is in fact nonlinear) is replaced by a sequence of linear boundary value problems (besides the first one)  with a singular linear operator called the exact functional discrete scheme (EFDS). A finite subsequence of $m$ terms, called truncated functional discrete scheme (TFDS), is the basis for our algorithm. The approach provides an
super-exponential convergence rate as $m \to \infty$. The eigenpairs can be
computed in parallel for all given indexes. The algorithm is  based on some recurrence procedures including the basic arithmetical operations with the coefficients of some expansions only. This is an exact symbolic algorithm (ESA) for $m=\infty$ and a truncated symbolic  algorithm (TSA) for a finite $m$. Numerical examples are presented to support the theory.
\end{abstract}
\date{July 10, 2006}
\end{topmatter}

\section{Introduction}

Using various definitions of fractional derivatives one can consider, e.g., differential operators of fractional
order, boundary and eigenvalue problems for these operators, as well as various approximation methods for
them (see e.g. \cite{mot, ford}).

The eigenvalue problem (EVP) is the problem of finding eigenvalues (frequencies) and eigenfunctions (vibration shapes) or so called eigenpairs  and
plays an important role in various applications concerned with
vibrations and wave processes \cite{af, pr}. Often it is needed to compute a
great number (hundreds of thousands) of eigenvalues and
eigenfunctions including eigenpairs with great indexes  (see e.g. \cite[p. 273]{pr}). A very efficient approach to solve such problems represent methods based on perturbation and  homotopy  ideas \cite{allg, arm}, such as the FD-method \cite{m1,m2, m3, dgm} (compare with the Adomian decomposition method \cite{a,r}). These approaches allow in a natural way the use of  computer algebra tools (see e.g. \cite{gg}).

In \cite{csw} a spectral approximation of Jacobi-type fractional differential equations (FDEs) is considered. Chen, Shen and Wang defined a
new class of generalized Jacobi functions (GJFs), which are the eigenfunctions of some fractional differential operator  and can serve as natural basis functions for properly designed spectral methods for FDEs. The spectral approximation results for these GJFs in weighted Sobolev spaces involving fractional derivatives are established and efficient GJF-Petrov-Galerkin methods for a class of prototypical
fractional initial value problems (FIVPs) and fractional boundary value problems
(FBVPs) of general order are constructed and analyzed. One of the important drawbacks  of these methods is the following: their accuracy decreaces with the growth of the eigenvalue index.

In order to find numerically the higher eigenvalues we propose a new approach  described below which we will  refer to  as the FD-method (following \cite{m1,m2, m3, dgm,gmr}). This  approach has been applied also to EVPs  with multiple eigenvalues \cite{gmr}. We would like to emphasize the following features of our approach, which differ from many other methods: 1) all eigenpairs can be computed in parallel, 2)The convergence rate increases as the index of the eigenpair increases, 3) the approximation called the truncated discrete scheme (TFDS) is derived from an exact discrete representation of the solution which we call the exact discrete scheme (EFDS), 4) the computational algorithm operates with the coefficients of some expansions of the initial data  and uses the basic arithmetical operations only so that the computer algebra tools are the natural medium of choice, 5) under some assumptions it is possible to pass to the limit at the truncated sum of the series and to obtain the exact solution in a closed form \cite{m4}. This is the reason, why we call our algorithm symbolic.

The article is organized as follows. First, we describe the FD-method for the fractional Jacobi-type differential operator. In Section \ref{s2} we prove the super-exponential convergence rate of our method for the case when the potential function is approximated by the constant zero. Section \ref{s3} deals with the reformulation of the algorithm in the case of a polynomial potential as a recursive procedure for the coefficients of the representation of corrections of the FD-method through the generalized Jacobi functions.  Although these corrections are the solutions of some BVPs with a singular differential operator, this procedure does not use any BVP-solver  but some recursions for the expansion coefficients with basic arithmetical operations only. We call such procedures symbolic algorithms. Note that similar algorithms for some special potential functions were proposed in \cite{dgm}. In Section \ref{s4}   numerical examples are given to support the theoretical results.

\section{Algorithm of the FD-method for a fractional Sturm-Liouville-type operator with an Jacobi-type weight function } \label{s1}

Let us define the left-side and right-side Riemann-Liouville fractional derivatives of order $s\in \left(0,1\right)$
\begin{equation} \label{1}
{}_{-1} D_{x}^{s} v(x)=\frac{1}{\Gamma (1-s)} \frac{d}{dx} \int _{-1}^{x}\frac{v(y)}{(x-y)^{s} } dy,\; \; \; \; \; x\in \left(-1,1\right),
\end{equation}
\begin{equation} \label{2}
{}_{x} D_{1}^{s} v(x)=-\frac{1}{\Gamma (1-s)} \frac{d}{dx} \int _{x}^{1}\frac{v(y)}{(y-x)^{s} } dy,\; \; \; \; \; x\in \left(-1,1\right),
\end{equation}
where $\Gamma \left(\alpha \right)$ is the gamma-function.
Let us consider the following eigenvalue problem of Sturm-Liouville-type
\begin{equation} \label{4}
\begin{split}
&{{}^{+} L_{\alpha ,\beta }^{2s} u(x)+\left(q(x)-\lambda \right)u(x)=0,\; \; x\in \left(-1,1\right),\; \; \, s\in \left(0,1\right),\; \; \;  \alpha \ge 0,\;\;\beta >-1,} \\
&{\left|u(-1)\right|<\infty ,\; \; \; \left|u(1)\right|<c \delta _{\alpha,0},\;\;\; c=\text{const} <\infty}
\end{split}
\end{equation}
for the fractional Jacobi-type differential operator
\begin{equation} \label{3}
{}^{+} L_{\alpha ,\beta }^{2s} u(x)=\omega ^{(\alpha ,-\beta )} {}_{-1} D_{x}^{s} \left\{\omega ^{(-\alpha +s,\beta +s)} {}_{x} D_{1}^{s} u(x)\right\},\; \; x\in \left(-1,1\right)
\end{equation}
with a piecewise smooth potential coefficient $q(x)$ and
the Jacobi-type weight function 
$$\omega ^{(\alpha ,\beta )} =\omega ^{(\alpha ,\beta )} (x)=(1-x)^{\alpha } (1+x)^{\beta }.$$
 In \eqref{4} $\delta _{n,k} $ denotes the Kronecker delta. A solution of the eigenvalue problem \eqref{4} consists of an eigenvalue $\lambda$ and of a corresponding eigenfunction $u(x)$, in other words, of the eigenpair $\lambda, u(x)$. If the set of eigenpairs is countable, we identify each pair by a natural number $n=1,2,...$. According to \cite{csw} (see Corollary 3.1) the Sturm-Liouville singular differential operator \eqref{3} is selfadjoint provided that $s\in \left(0,1\right)$. For $s-1<\alpha<0$ the second boundary condition in \eqref{4} is given as follows:
\begin{equation} \label{3-0}
\left|\omega ^{(-\alpha +s,\beta +s)} {}_{x} D_{1}^{s} \left\{(1-x)^{\alpha } P_{n}^{(\alpha ,\beta )} (x) \right\}\right|_{x=1}<\infty.
\end{equation}

We apply to problem (\ref{4}) the FD-method (see e.g. \cite{dgm} for details) which is a combination of perturbation of the original differential operator by a parameter dependent operator (embedding) and then a ``trip'' along this parameter from a ``simple'' problem to the original one (homotopy). The homotopy idea was exploited in various ways e.g. in \cite{allg, arm}. The FD-method consists of the following steps: 1) the function $q(x)$ (the potential) is approximated by a piece-wise constant function $\overline{q}(x)$; 2) problem (\ref{4}) is embedded into the parametric family of problems with the potential  $\overline{q}(x)+t(q(x)-\overline{q}(x))$ (the perturbation of the operator); 3) the solution of the perturbed family is represented as a power series in $t$ with the coefficients which are solutions of a recursive sequence of problems with the potential $\overline{q}(x)$; 4) by setting $t=1$ one obtains the series representations of the solution of the original problem. Below we consider the simplest case of the approximation of $q(x)$ by the constant $\overline{q}(x)\equiv 0$. Note that in this case the idea of the FD-method is close to the idea of the Adomian method \cite{a,r}.

The exact solution of the eigenvalue problem is then represented by the series
\begin{equation} \label{5-0}
{u}_{n}(x)=\underset{j=0}{\overset{\infty}{\sum}} u_{n}^{(j)}(x), \; \; \; {\lambda}_{n}=\underset{j=0}{\overset{\infty}{\sum}} \lambda_{n}^{(j)}
\end{equation}
provided that these series converge.
The sufficient conditions for the convergence of the series \eqref{5-0} will be presented later in Section~\ref{s2} (see  \eqref{34}).
The approximate solution to problem \eqref{4} is represented by a pair of corresponding truncated series
\begin{equation} \label{5}
\overset{m}{u}_{n}(x)=\underset{j=0}{\overset{m}{\sum}} u_{n}^{(j)}(x), \; \; \; \overset{m}{\lambda}_{n}=\underset{j=0}{\overset{m}{\sum}} \lambda_{n}^{(j)}
\end{equation}
which is called an approximation of rank $m$. The summands of series (\ref{5-0}), (\ref{5}) are the solutions of the following recursive sequence of problems:
\begin{equation} \label{5-1}
\begin{split}
 &{{}^{+} L_{\alpha ,\beta }^{2s} u_{n}^{(j+1)} (x)-\lambda _{n}^{(0)} u_{n}^{(j+1)} (x)=F_{n}^{(j+1)} (x),\; \; x\in \left(-1,1\right),} \\
 &{s\in \left(0,1\right),\; \; \; \alpha \ge 0,\;\;\;\beta >-1,} \\
 &{\left|u_{n}^{(j+1)} (-1)\right|<\infty ,\; \; \; \; \; \left|u_{n}^{(j+1)} (1)\right| \le c \delta _{\alpha,0},\;\;\; c=const<\infty,}
 \end{split}
\end{equation}
where
\begin{equation} \label{5-2}
F_{n}^{(j+1)} (x)=\sum _{p=0}^{j}\lambda _{n}^{(j+1-p)}  u_{n}^{(p)} (x)-q(x)u_{n}^{(j)} (x),\; \; \; j=0,1,2,...,\; \; \; n=0,1,2,...
\end{equation}
and $\delta _{\alpha, \beta}$ is the Kronecker delta.

For $s-1<\alpha<0$ the second boundary condition in \eqref{5-1} is given as follows:
\begin{equation} \label{5-3}
\left|\omega ^{(-\alpha +s,\beta +s)} {}_{x} D_{1}^{s} \left\{u_{n}^{(j+1)} (x) \right\}\right|_{x=1}<\infty.
\end{equation}
The initial approximation $u_{n}^{(0)} (x),$ $\lambda _{n}^{(0)} $ is the solution of the so called base problem, that is,
\begin{equation} \label{6}
\begin{split}
&{{}^{+} L_{\alpha ,\beta }^{2s} u_{n}^{(0)} (x)-\lambda _{n}^{(0)} u_{n}^{(0)} (x)=0,\; \; x\in \left(-1,1\right),\, \; \; \; s\in \left(0,1\right),\; \;\;  \alpha \ge 0,\; \; \;\beta >-1} \\
&{\left|u_{n}^{(0)} (-1)\right|<\infty ,\; \; \; \; \; \left|u_{n}^{(0)} (1)\right| \le c \delta _{\alpha,0},\;\;\; c=const<\infty}
\end{split}
\end{equation}
and for $s-1<\alpha<0$ the second boundary condition in \eqref{6} is given by
\begin{equation} \label{6-1}
\left|\omega ^{(-\alpha +s,\beta +s)} {}_{x} D_{1}^{s} \left\{u_{n}^{(0)} (x) \right\}\right|_{x=1}<\infty.
\end{equation}
We call the representation of the exact solution by the series (\ref{5-0}), i.e., by the sequences $\{u_{n}^{(j)}, j=0,1,...\}$, and we call $\{\lambda_{n}^{(j)},  j=0,1,...\}$  we call the exact discrete scheme (EDisS) and the approximate representation by the corresponding truncated series the truncated discrete scheme (TDisS).
The idea in broad sense is related to the idea of exact and truncated difference schemes (EDS and TDS) from \cite{ghmk} but the algorithms are rather different.

The solution of (\ref{6}) is the generalized Jacobi function
\begin{equation} \label{7}
{}^{+} J_{n}^{(-\alpha ,\beta )} (x)=(1-x)^{\alpha } P_{n}^{(\alpha ,\beta )} (x),
\end{equation}
where $P_{n}^{(\alpha ,\beta )} (x)$  is the classical Jacobi polynomial. This is a generalization of the eigenfunctions of the classical Jacobi differential operator (see (\ref{10}) and Theorem 3.2 in \cite{csw}). Since the operator in (\ref{5-1}) is singular these problems are solvable under the solvability condition $$\int _{-1}^{1} F_{n}^{(j+1)} (x)u_{n}^{(0)} (x)\omega^{(-\alpha ,\beta )} (x)dx=0.$$
The solution of (\ref{6}) can be found up to a constant
whose value is determined by the  orthogonality condition
\begin{equation} \label{7-1}
\int _{-1}^{1} u_{n}^{(0)} (x)u_{n}^{(j+1)} (x)\omega ^{(-\alpha ,\beta )} (x)dx=0.
\end{equation}
The polynomials $P_{n}^{(\alpha ,\beta )}(x)$ are orthogonal with the Jacobi weight function  $\omega ^{(\alpha ,\beta )}(x)$ (see \cite{gr}), i.e.
\begin{equation} \label{8}
\int _{-1}^{1}P_{n}^{(\alpha ,\beta )} (x)P_{k}^{(\alpha ,\beta )} (x) \omega ^{(\alpha ,\beta )} (x)dx=\gamma _{n}^{(\alpha ,\beta )} \delta _{n,k} ,
\end{equation}
where
\begin{equation} \label{9}
\gamma _{n}^{(\alpha ,\beta )} =\frac{2^{\alpha +\beta +1} \Gamma (n+\alpha +1)\Gamma (n+\beta +1)}{(2n+\alpha +\beta +1)n!\Gamma (n+\alpha +\beta +1)}.
\end{equation}
The normalized solution of problem (\ref{6}) is
\begin{equation} \label{10}
\begin{split}
&{u_{n}^{(0)} (x)=\left(\gamma _{n}^{(\alpha ,\beta )} \right)^{-1/2} \cdot {}^{+} J_{n}^{(-\alpha ,\beta )} (x),} \\
&{\lambda _{n}^{(0)} =\frac{\Gamma (n+\alpha +1)\Gamma (n+\beta +s+1)}{\Gamma (n+\alpha -s+1)\Gamma (n+\beta +1)} ,\; \; \; n=0,1,2,...}
\end{split}
\end{equation}
Using the solvability condition and the orthogonality property (\ref{7-1}) we obtain from (\ref{5-1}), (\ref{5-2}) the following formula for the eigenvalue corrections
\begin{equation} \label{13}
\lambda _{n}^{(j+1)} =\int _{-1}^{1} q(x)u_{n}^{(0)} (x)u_{n}^{(j)} (x)\omega ^{(-\alpha ,\beta )} (x)dx.
\end{equation}
We consider the solution $u_{n}^{(j+1)} (x)$ of problem (\ref{5-1}), (\ref{5-2}) satisfying the orthogonality condition (\ref{7-1}) as an element of the Hilbert space $L_{w}^{2} \left[-1,1\right]$ with the scalar product $(u,v)=\int_{-1}^1\omega^{(-\alpha ,\beta )}(x)u(x)v(x)dx$.
Now, using (\ref{7-1})
we obtain the following representation for the eigenfunction corrections
\begin{equation} \label{15}
u_{n}^{(j+1)} (x)=\underset{p=0,p\ne n}{\overset{\infty }{\sum}}\int _{-1}^{1} F_{n}^{\left(j+1\right)} \left(\xi \right)u_{p}^{(0)} \left(\xi \right)\omega ^{(-\alpha ,\beta )} \left(\xi \right)d\xi  \frac{u_{p}^{(0)} (x)}{\lambda _{n}^{(0)} -\lambda _{p}^{(0)} } .
\end{equation}

\begin{remark}
Note that (\ref{5-0}) can be considered as an exact discrete scheme for the boundary value problem (\ref{4}) where we have switched from a nonlinear problem to the sequence $u_{n}^{(j)}(x)$ of solutions of linear problems.
The truncated series (\ref{5}) can be viewed as an algorithmically realizable truncated discrete scheme.
\end{remark}

\section{Convergence of the FD-method with $\overline{q}(x)\equiv 0$} \label{s2}

In this section we find sufficient convergence conditions of the FD-method as well as the accuracy estimates.

From (\ref{15}) we obtain
\begin{equation} \label{16}
\begin{split}
 &{\left\| u_{n}^{(j+1)} \right\| \; =\left(\underset{p=0,p\ne n}{\overset{\infty}{\sum}}\left(\int _{-1}^{1} F_{n}^{\left(j+1\right)} \left(\xi \right)u_{p}^{(0)} \left(\xi \right)\omega^{(-\alpha ,\beta )} \left(\xi \right)d\xi \left(\lambda _{n}^{(0)} -\lambda _{p}^{(0)} \right)^{-1} \right)^{2}  \right)^{1/2} } \\
 &{\le \max \left\{\left(\lambda _{n}^{(0)} -\lambda _{n-1}^{(0)} \right)^{-1} ,\left(\lambda _{n+1}^{(0)} -\lambda _{n}^{(0)} \right)^{-1} \right\}\left\| F_{n}^{\left(j+1\right)} \right\| =M_{n} \left\| F_{n}^{\left(j+1\right)} \right\| ,}
 \end{split}
\end{equation}
where
\begin{equation} \label{17}
\left\| u_{n}^{(j+1)} \right\|=\left\| u_{n}^{(j+1)} \right\|_{L_{\omega}^2} =\left(\int _{-1}^{1} \left(u_{n}^{(j+1)} \left(x\right)\right)^{2} \omega^{(-\alpha ,\beta )} \left(x\right)dx\right)^{{1\mathord{\left/ {\vphantom {1 2}} \right. \kern-\nulldelimiterspace} 2} } ,       \left\| u_{n}^{(0)} \right\| =1,
\end{equation}
and
\begin{equation} \label{18}
\begin{split}
&{M_{n} =\max \left\{\frac{1}{\lambda _{n}^{(0)} -\lambda _{n-1}^{(0)} } ,\frac{1}{\lambda _{n+1}^{(0)} -\lambda _{n}^{(0)} } \right\}} \\
&{=\frac{\Gamma \left(n+\beta +1\right)\Gamma \left(n+\alpha +1-s\right)}{s\left(2n+\beta +\alpha \right)\Gamma \left(n+\alpha \right)\Gamma \left(n+\beta +s\right)}}\\
&{\times \max \left\{1,\frac{\left(n+\beta +1\right)\left(2n+\beta +\alpha \right)\left(n+\alpha +1-s\right)}{\left(2n+\beta +\alpha +2\right)\left(n+\alpha \right)\left(n+\beta +s\right)} \right\}.}
\end{split}
\end{equation}
Here we have $M_{n}=1$ if coincidentally $s=1/2$ and $\alpha =\beta$ or $s=1/2$ and $\alpha =\beta +1$.

The asymptotic formula from \cite{te}, namely, 
\begin{equation} \label{19}
\frac{\Gamma \left(z+\gamma \right)}{\Gamma \left(z+\eta \right)} =z^{\gamma -\eta } \left[1+\frac{\left(\gamma -\eta \right)\left(\gamma +\eta -1\right)}{2z} +o\left(\left|z\right|^{-2} \right)\right]\,\,\, \text{as}\,\,\, z \to \infty
\end{equation}
together with (\ref{18}) yields
\begin{equation}\label{27}
\begin{split}
&\underset{n \to \infty}{\lim}M_{n}=\frac{1}{2s}\cdot \left\{
                                                        \begin{array}{ll}
                                                          \infty, & \; s \in (0,1/2),  \\
                                                          1, & \; s=1/2, \\
                                                          0, & \; s \in (1/2,1).
                                                        \end{array}
                                                      \right.
\end{split}
\end{equation}

From (\ref{16}) and (\ref{13}) we obtain
\begin{equation} \label{28}
\left\| u_{n}^{(j+1)} \right\| \le M_{n} \left\| q\right\| _{\infty } \left\{\sum _{l=1}^{j} \left\| u_{n}^{(j-l)} \right\| \left\| u_{n}^{(l)} \right\| +\left\| u_{n}^{(j)} \right\| \right\},\; \; \; \; {\kern 1pt} {\kern 1pt} \left|\lambda _{n}^{(j+1)} \right|\le \left\| q\right\| _{\infty } \left\| u_{n}^{(j)} \right\| ,
\end{equation}
where $\left\| q\right\| _{\infty } =\mathop{\max }\limits_{x\in \left[-1,1\right]} \left|q(x)\right|$. The estimates (\ref{28}) lead to
\begin{equation} \label{29}
\left\| u_{n}^{(j+1)} \right\| \le M_{n} \left\| q\right\| _{\infty } \sum _{l=0}^{j} \left\| u_{n}^{(j-l)} \right\| \left\| u_{n}^{(l)} \right\| .
\end{equation}
Substituting here $U_{j} =\left(\left\| q\right\| _{\infty } M_{n} \right)^{-j} \left\| u_{n}^{\left(j\right)} \right\|$ and $U_{0} =\left\| u_{n}^{\left(0\right)} \right\| =1$ and replacing the new variables by the majorant variables subject to $U_{j} \le \overline{U}_{j} $ and $\overline{U}_{0} =U_{0} =1$, we come to the majorant equation
\begin{equation} \label{30}
\overline{U}_{j+1} =\sum _{l=0}^{j} \; \overline{U}_{j-l} \overline{U}_{l}.
\end{equation}
The solution of this convolution type equation is (see e.g. \cite[p.~159-161,210]{v} and \cite{rnd})
\begin{equation} \label{31}
\overline{U}_{j+1} =\frac{\left(2j+2\right)!}{\left(j+1\right)!\left(j+2\right)!} =4^{j+1} 2\frac{\left(2j+1\right)!!}{\left(2j+4\right)!!} .
\end{equation}
Returning to the old variables we obtain the following estimate for the solution of (\ref{29})
\begin{equation} \label{32}
\left\| u_{n}^{\left(j+1\right)} \right\| \le \left(4\left\| q\right\| _{\infty } M_{n} \right)^{j+1} 2\frac{\left(2j+1\right)!!}{\left(2j+4\right)!!} \le \frac{\left(4\left\| q\right\| _{\infty } M_{n} \right)^{j+1} }{\left(j+2\right)\sqrt{\pi \left(j+1\right)} } ,
\end{equation}
and then from (\ref{28}) the next estimate for the eigenvalue corrections
\begin{equation} \label{33}
\left|\lambda _{n}^{(j+1)} \right|\le \left\| q\right\| _{\infty }\left(4\left\| q\right\| _{\infty } M_{n} \right)^{j} 2\frac{\left(2j-1\right)!!}{\left(2j+2\right)!!} \le \left\| q\right\| _{\infty }\frac{\left(4\left\| q\right\| _{\infty } M_{n} \right)^{j} }{\left(j+1\right)\sqrt{\pi j} } .
\end{equation}
The last part of inequalities (\ref{32}) and (\ref{33}) was obtained using the reflections like those from the proof of the Wallis formula (see e.g. \cite[p.~344]{f}).

From (\ref{32}) and (\ref{33}) follows the next assertion.

\begin{theorem}
Under the condition
\begin{equation} \label{34}
r_{n} =4\left\| q\right\| _{\infty } M_{n} <1,  n=0,1,2,...
\end{equation}
the FD-method converges super-exponentially and his accuracy is characterized by the estimates
\begin{equation} \label{35}
\left|\lambda _{n} -\overset{m}{\lambda}_{n} \right|\le 2 \left\| q\right\| _{\infty } \frac{\left(r_{n} \right)^{m} }{1-r_{n} } \frac{\left(2m-1\right)!!}{\left(2m+2\right)!!} \le \left\| q\right\| _{\infty } \frac{\left(r_{n} \right)^{m} }{1-r_{n} } \frac{1}{\left(m+1\right)\sqrt{\pi m} } 
\end{equation}
and
\begin{equation} \label{36}
\left\| u_{n} -\overset{m}{u}_{n} \right\| \le 2\frac{\left(r_{n} \right)^{m+1} }{1-r_{n} } \frac{\left(2m+1\right)!!}{\left(2m+4\right)!!} \le \frac{\left(r_{n} \right)^{m+1} }{1-r_{n} } \frac{1}{\left(m+2\right)\sqrt{\pi \left(m+1\right)} } .
\end{equation}
\end{theorem}

\begin{remark}
The relation (\ref{27}) shows that for a fixed $s \in (1/2,1)$ there exists some $n_{0}$ such that for all $n\ge n_{0} $ condition  (\ref{34}) of the theorem is fulfilled and the  FD-method converges super-exponentially. The error estimates (\ref{35}) and (\ref{36}) show that the accuracy of our method increases as we increase the eigenvalue number. This remarkable property can be lost for problem (\ref{4}) with a fixed $0<s<1/2$  or  if coincidentally $s=1/2$ and $\alpha =\beta$ or $s=1/2$ and $\alpha =\beta +1$. For these cases condition (\ref{34}) for a given $n$ can be fulfilled for $\left\| q\right\| _{\infty } $ small enough. If this is not the case then the trivial variant of the FD-method with  $\overline{q}(x)\equiv 0$ is divergent and one should apply the variant with a piecewise constant $\overline{q}(x)$ (see  \cite{m1, m2}).
\end{remark}

\section{Symbolic algorithm of the FD-method with $\overline{q}(x)\equiv 0$} \label{s3}

  A new algorithmic realization of the FD-method for problem  (\ref{4}) with $\bar{q}(x)\equiv 0, \alpha =\beta =0$ and with the potential $q(x)=x^{2} $ was proposed in \cite{dgm}. It was shown that the corrections to eigenfunctions of the FD-method are linear combinations of Legendre polynomials with a number of summands depending on the degree of the potential polynomial and on the correction number $j$. The coefficients of these linear combinations can be represented recursively through the corresponding coefficients  computed at previous steps  using the basic arithmetical operations only. The approximations of the eigenvalues are represented through these coefficients too.  Thus, the method uses these arithmetical operations and does not involve  solutions of any supplementary BVPs and computation of any integrals unlike the traditional implementations. In this sense our algorithm is a symbolic one since it operates with the expansion coefficients only.

Further we extend these results and describe such algorithm for problem (\ref{4}) with the potential
\begin{equation} \label{37}
q(x)=\sum _{l=0}^{r}c_{l} x^{l}\,\,\,\,\,   (c_{r} \ne 0),
\end{equation}
where $c_{l}$, $l=0,1,...,r$ are real constant coefficients. In this case the FD-method is exactly realizable in the sense that the corrections to the eigenfunctions can be explicitly represented as linear combinations of the generalized Jacobi functions ${}^{+} J_{n}^{(-\alpha ,\beta )} (x)$ with coefficients which can be explicitly represented as polynomials of $c_{l}$, $l=0,1,...,r$ with rational coefficients.

The solution $u_{n}^{(0)} (x)$ of the base problem (\ref{10}) (up to an arbitrary multiplicative constant) can be represented as  $u_{n}^{(0)} (x)=a_{n}^{(0)} {}^{+} J_{n}^{(-\alpha ,\beta )} (x)$, $a_{n}^{(0)} \ne 0$ (see Lemma 4.1 below). Note that we can renounce here the normalizing of the eigenfunctions of the base problem which we have used above to prove the convergence of the FD-method.  For the orthonormal eigenfunctions $\left\{u_{n}^{(0)} (x)\right\}_{n=1}^{\infty } $ in accordance with (\ref{10}) we have $a_{n}^{(0)} =\left(\gamma _{n}^{(\alpha ,\beta )} \right)^{-1/2} $.

Further we use the following recurrence relation for the Jacobi polynomials (see, e.g., \cite[Section 10.8]{be})
\begin{equation} \label{38}
x^{r} P_{n}^{(\alpha ,\beta )} (x)=\sum _{k=\max (n-r,0)}^{n+r}b_{n+r+1-k,n,r}  P_{k}^{(\alpha ,\beta )} (x),\,\,\,\,\,\,\, r,n=0,1,2,...,
\end{equation}
where
\[b_{1,n,0} =1,\,\,\,\,\,\,\,    b_{1,n,1} =\frac{2(n+1)(n+\alpha +\beta +1)}{(2n+\alpha +\beta +1)(2n+\alpha +\beta +2)} ,\]
\[b_{2,n,1} =\frac{\beta ^{2} -\alpha ^{2} }{(2n+\alpha +\beta +1)(2n+\alpha +\beta +2)} ,\]
\[b_{3,n,1} =\frac{2(n+\alpha )(n+\beta )}{(2n+\alpha +\beta +1)(2n+\alpha +\beta )} ;\]
\[b_{n+r+1-k,n,r} =\sum _{t=\max (0,k-1,n-r+1)}^{\min (k+1,n+r-1)}b_{n+r-t,n,r-1} b_{t-k+2,t,1}  ,\]
with $k=\max (n-r,0),\max (n-r,0)+1,...,n+r$ and $r=2,3,...$.
The $(j+1)$-th step of the FD-method consists of solving the BVP (\ref{5-1}), (\ref{5-2}) with
\begin{equation} \label{39}
F_{n}^{(j+1)} (x)=\sum _{p=0}^{j}\lambda _{n}^{(j+1-p)}  u_{n}^{(p)} (x)-\sum _{l=0}^{r}c_{l} x^{l}  u_{n}^{(j)} (x).
\end{equation}
Using the solution of the base problem (\ref{10}) and considering the properties of the classical Jacobi polynomials $P_{n}^{(\alpha ,\beta )} (x)$, of the generalized Jacobi functions ${}^{+} J_{k}^{(-\alpha ,\beta )} (x)$, the properties of the fractional Sturm-Liouville type operator ${}^{+} L_{\alpha ,\beta }^{2s} u(x)$ as well as the relations (\ref{38}) and (\ref{39}) we come to the following assertion.

\begin{lemma}
The solution of problem (\ref{5-1}) with the right-hand sides (\ref{39}) can be represented by
\begin{equation} \label{40}
u_{n}^{(j+1)} (x)=\underset{k=\max (0,\; n-r(j+1))}{\overset{n+r(j+1)}{\sum}}a_{k}^{(j+1)}  \; {}^{+} J_{k}^{(-\alpha ,\beta )} (x),\, \, \, \; \; n,j=0,1,2,...,
\end{equation}
where $a_{n}^{(j+1)} =0$, $j=0,1,2,...$, and $a_{n}^{(0)}\neq 0$.
\end{lemma}

Let us substitute (\ref{40}) into (\ref{39}) and use (\ref{38}) with the aim to represent (\ref{39}) through the generalized Jacobi functions  ${}^{+} J_{k}^{(-\alpha ,\beta )} (x)$ only. Then, by changing the summation order we arrive at
\begin{equation} \label{41}
\begin{split}
 &{F_{n}^{(j+1)} (x)=\left(\lambda _{n}^{(j+1)} a_{n}^{(0)} -\underset{t=\max (0,n-r)}{\overset{n+r}{\sum}}a_{t}^{(j)} \underset{l=\max (0,n-t,t-n)}{\overset{r}{\sum}} c_{l} b_{t+l+1-n,t,l}   \right) \; {}^{+} J_{n}^{(-\alpha ,\beta )} (x)} \\
 &{+\underset{m=\max (0,n-rj), m\ne n}{\overset{n+rj}{\sum}}{}^{+} J_{m}^{(-\alpha ,\beta )} (x)\underset{p=\left[\frac{|m-n|+r-1}{r} \right]}{\overset{\min (j,n+rj)}{\sum}}\lambda _{n}^{(j+1-p)}  a_{m}^{(p)}  } \\
 &{-\underset{m=\max (0,\; n-r(j+1)), m\ne n}{\overset{n+r(j+1)}{\sum}}{}^{+} J_{m}^{(-\alpha ,\beta )} (x)\underset{l=0}{\overset{r}{\sum}}c_{l}  \underset{k=\max (0,n-rj,m-l)}{\overset{\min (n+rj,m+l)}{\sum}}a_{k}^{(j)} b_{k+l+1-m,k,l}  ,}
 \end{split}
\end{equation}
where $\left[\sigma \right]$ is the entire part of $\sigma $. Substituting (\ref{40}) into (\ref{7-1}) and comparing the coefficients in the front of the generalized Jacobi functions ${}^{+} J_{k}^{(-\alpha ,\beta )}(x),$ where $k$ is an integer in
$\left[\max (0,\; n-r(j+1)),n+r(j+1)\right],$
 we obtain the following formulas for the eigenvalue corrections
\begin{equation} \label{42}
\lambda _{n}^{(j+1)} =\left(a_{n}^{(0)} \right)^{-1} \; \sum _{t=\max (0,n-r)}^{n+r}a_{t}^{(j)} \sum _{l=\max (0,n-t,t-n)}^{r}c_{l} b_{t+l+1-n,t,l}   ,\; \; j=0,1,2,...,
\end{equation}
\[a_{n}^{(j)} =0,\; j=1,2,...,\; \]
\begin{equation} \label{43}
\lambda _{n}^{(0)} =\frac{\Gamma (n+\alpha +1)\Gamma (n+\beta +s+1)}{\Gamma (n+\alpha -s+1)\Gamma (n+\beta +1)}
\end{equation}
as well as the following recursive representations for the coefficients in  (\ref{40}):
\begin{equation} \label{44}
\begin{split}
&a_{m}^{(j+1)} =-\left(\lambda _{m}^{(0)} -\lambda _{n}^{(0)} \right)^{-1} \underset{l=0}{\overset{r}{\sum}}c_{l} \; \underset{k=\max (0,n-rj,m-l)}{\overset{\min (n+rj,m+l)}{\sum}}a_{k}^{(j)} b_{k+l+1-m,k,l},
\end{split}
\end{equation}
with $m=n+rj+1,n+rj+2,...,n+r(j+1),\; \; m\ne n,\; \; j=0,1,2,...$,
\begin{equation} \label{45}
\begin{split}
&a_{m}^{(j+1)} =\left(\lambda _{m}^{(0)} -\lambda _{n}^{(0)} \right)^{-1} \left(\underset{p=\left[\frac{|m-n|+r-1}{r} \right]}{\overset{\min (j,n+rj)}{\sum}}\lambda _{n}^{(j+1-p)}  a_{m}^{(p)} -\underset{l=0}{\overset{r}{\sum}}c_{l}  \sum _{k=\max (0,n-rj,m-l)}^{\min (n+rj,m+l)}a_{k}^{(j)} b_{k+l+1-m,k,l}  \right),
\end{split}
\end{equation}
with $m=\max (0,n-rj),\max (0,n-rj)+1,...,n+rj,\; \; m\ne n,\; \; j=0,1,2,...$,
\begin{equation} \label{46}
\begin{split}
&a_{m}^{(j+1)} =-\left(\lambda _{m}^{(0)} -\lambda _{n}^{(0)} \right)^{-1} \sum _{l=0}^{r}c_{l}  \sum _{k=\max (0,n-rj,m-l)}^{\min (n+rj,m+l)}a_{k}^{(j)} b_{k+l+1-m,k,l},
\end{split}
\end{equation}
with $m=\max (0,n-r(j+1)),\max (0,n-r(j+1))+1,...,\max (0,n-rj-1),\; \;m\ne n,\; \; j=0,1,2,...$,
and
\begin{equation} \label{47}
a_{n}^{(j+1)} =0,\;\;\;\; j=0,1,2,...,\;\;\;\;a_{n}^{(0)}\neq 0.
\end{equation}
Formulas (\ref{40})-(\ref{47}) represent the symbolic algorithm of the FD-method for problem  (\ref{4}) with the polynomial potential (\ref{37}).

\section{Numerical examples}\label{s4}

\begin{example} \label{ex_r1}
 We consider the eigenvalue problem (\ref{4}) with  $\alpha =1/2 ,\beta =0,$ $s=3/4$, and with  potential (\ref{37}),  where  $r=3$, $c_{3} =1/4$ and $c_{l} =0$, $l=0,1,2$. In this case the sufficient convergence condition  $r_{n}<1$ (see (\ref{34})) is fulfilled for $n\ge 1$ (see Table~1). We have $\left\| q\right\| _{\infty } =1/4$ and in \eqref{18}
\[M_{n} =\frac{8}{3} \frac{\Gamma \left(n+1\right)}{\left(4n+1\right)\Gamma \left(n+1/2 \right)}.\]
The computer algebra system Maple~17 was used, where the corrections \eqref{40} and \eqref{42} to the eigenpairs were computed exactly as mathematical expressions, i.e., we had no rounding errors. Below we give the eigenvalue corrections \eqref{42} and the coefficients of (\ref{40}) for some first steps of the FD-method:
\[\lambda _{0}^{(0)} =\frac{3}{8} \sqrt{\pi } ,\; \; \; \; u_{0}^{(0)} (x)=a_{0}^{(0)} {}^{+} J_{0}^{(-1/2 ,0 )} (x),\]
\[\lambda _{0}^{(1)} =-\frac{13c_{3} }{105} , \lambda _{0}^{(2)} =-\frac{{\rm 201134942464}}{{\rm 1943987920875}} \frac{\left(c_{3} \right)^{2} }{\sqrt{\pi } } ,\;  \lambda _{0}^{(3)} =-\frac{{\rm 274356801766461046784}}{{\rm 81295088830639587028125}} \frac{\left(c_{3} \right)^{3} }{\pi } ,\; \]
\[\lambda _{0}^{(4)} \approx {\rm 0.001622...}\cdot \frac{\left(c_{3} \right)^{4} }{\pi ^{{3\mathord{\left/ {\vphantom {3 2}} \right. \kern-\nulldelimiterspace} 2} } } ,\; \dots , \lambda _{0}^{(10)} \approx {\rm 0.000005157...}\cdot \frac{\left(c_{3} \right)^{10} }{\pi ^{{9\mathord{\left/ {\vphantom {9 2}} \right. \kern-\nulldelimiterspace} 2} } }, \]
\[a_{0}^{(0)} =1,\,\,\,a_{1}^{(1)} =-\frac{1216}{2475} \frac{a_{0}^{(0)} c_{3} }{\sqrt{\pi } } ,\,\,\, a_{2}^{(1)} =\frac{2048}{38493} \frac{a_{0}^{(0)} c_{3} }{\sqrt{\pi } } ,\,\,\, a_{3}^{(1)} =-\frac{16384}{204633} \frac{a_{0}^{(0)} c_{3} }{\sqrt{\pi } } ,\]
\[a_{1}^{(2)} =\frac{{\rm 290443255808}}{{\rm 70816702831875}} \frac{a_{0}^{(0)} \left(c_{3} \right)^{2} }{\pi } , a_{2}^{(2)} =\frac{{\rm 193813841149952}}{{\rm 1609728034189275}} \frac{a_{0}^{(0)} \left(c_{3} \right)^{2} }{\pi } ,\]
\[a_{3}^{(2)} =-\frac{{\rm 48103527022592}}{{\rm 8557490370203775}} \frac{a_{0}^{(0)} \left(c_{3} \right)^{2} }{\pi } , a_{4}^{(2)} =\frac{{\rm 340833471561728}}{{\rm 13515289050237975}} \frac{a_{0}^{(0)} \left(c_{3} \right)^{2} }{\pi } ,\]
\[a_{5}^{(2)} =-\frac{{\rm 336081190912}}{{\rm 249927686251425}} \frac{a_{0}^{(0)} \left(c_{3} \right)^{2} }{\pi } , a_{6}^{(2)} =\frac{{\rm 68719476736}}{{\rm 48376171671975}} \frac{a_{0}^{(0)} \left(c_{3} \right)^{2} }{\pi } .\]
The approximations $\overset{m}{\lambda }_{n} $ of rank $m$  to the eigenvalues $\lambda _{n} $ for
$n=0,1,2,3,4,10$ using the FD-method of rank $m=20$ are given in Table 1.
The behavior of the corrections $\lambda _{n}^{(m)} ,\; u_{n}^{(m)} (x)$,  $m=0,1,...,10,20$ for the eigenpairs $\lambda _{n} ,$ $u_{n} (x),$ $n=0,10,$ is illustrated in  Table~2, which contains the corrections to these eigenvalues as well as the norms \eqref{17} of the eigenfunction corrections $|| u_{n}^{(m)} ||$.  One can observe  that the convergence of the method improves when the eigenpair number  $n$ increases (see Remark 3.1). Table 2 shows that the method converges for $n=0$ too, i.e., the conditions of Theorem 3.1 are rough and can be improved.
\begin{table}
\begin{center}
\begin{tabular}{|c|c|c|} \hline
$n$ & $r_{n}$ & $\overset{20}{\lambda }_{n} $ \\ \hline \hline
0 & 1.505 & 0.630053891717269391713596782178 \\ \hline
1 & 0.602 & 2.30514376782605437183729542346 \\ \hline
2 & 0.446 & 4.56663095405447807274309207123 \\ \hline
3 & 0.370 & 7.26852041545097231248033107978 \\ \hline
4 & 0.324 & 10.3587519350445765031940934458 \\ \hline
10 & 0.208 & 35.2508805155975526382623041970 \\ \hline
\end{tabular}
\caption{Approximations $\overset{20}{\lambda }_{n}$ to the eigenvalues $\lambda _{n}$ and the values of the ratio  $r_{n}$ of the geometric progression for $n=0,1,2,3,4,10$.}
\end{center}
\end{table}

\begin{table}
\begin{center}
\begin{tabular}{|c|c|c|c|c|} \hline
$m$  & $\lambda _{0}^{(m)} $ & $|| u_{0}^{(m)} ||  $ & $\lambda _{10}^{(m)} $ & $|| u_{10}^{(m)} ||$ \\ \hline \hline
0 & 0.6647 &  1.373 & 35.25 &  0.3627 \\ \hline
1 & $-0.03095$ & 0.06253 & $-0.0002053$ &  0.009806 \\ \hline
2 & $-0.003648$ &  0.001742 & 0.00009424 &  0.00007673 \\ \hline
3 & $-0.00001679$ &  0.00004165 &$ -1.207 \cdot {10}^{-8}$ &  0.000002101 \\ \hline
4 & 0.000001138  & 0.000002067 & $6.993 \cdot {10}^{-10}$  & $2.059 \cdot {10}^{-8}$ \\ \hline
5 & $9.723 \cdot {10}^{-8}$& $4.556 \cdot {10}^{-8}$& $-7.034 \cdot {10}^{-13}$ & $5.656 \cdot {10}^{-10}$\\ \hline
6 &$ -1.319 \cdot {10}^{-10}$ & $5.386 \cdot {10}^{-9}$ & $1.597 \cdot {10}^{-14}$ & $5.657 \cdot {10}^{-12}$ \\ \hline
7 &  $-2.950 \cdot {10}^{-10}$ & $ 1.856 \cdot {10}^{-10} $& $-4.932 \cdot {10}^{-17}$ & $1.553 \cdot {10}^{-13}$  \\ \hline
8 &$ -4.443 \cdot {10}^{-12}$ &  $1.438 \cdot {10}^{-11}$ & $7.277 \cdot {10}^{-19}$ & $1.558 \cdot {10}^{-15}$ \\ \hline
9 &  $7.229 \cdot {10}^{-13}$  & $7.175 \cdot {10}^{-13}$ & $-6.316 \cdot {10}^{-21}$ & $4.277 \cdot {10}^{-17}$ \\ \hline
10 & $2.848 \cdot {10}^{-14}$  & $3.451 \cdot {10}^{-14}$ &$ 2.910 \cdot {10}^{-23}$ & $4.293 \cdot {10}^{-19}$ \\ \hline
20 &$ -3.979 \cdot {10}^{-27}$ & $2.621 \cdot {10}^{-26}$ & $-3.694 \cdot {10}^{-45}$ & $6.793 \cdot {10}^{-35}$ \\ \hline
\end{tabular}

\caption{Corrections $\lambda^{(m)}_{n}, n=0,10$, and the norms of corrections to the corresponding  eigenfunctions $u^{(m)}_{n}, n=0,10$, for the FD-method of the ranks $m=0,1,...10, 20$. }
\end{center}
\end{table}

\end{example}

\begin{example} \label{ex_r2}
 Let us consider the eigenvalue problem (\ref{4}) in the case $s-1<\alpha<0$ with the second boundary condition \eqref{3-0}, namely with $\alpha=-1/8$, $\beta=-1/2$, $s =3/4$ and with  potential \eqref{37},  where  $r=3$, $c_{l} =1/12$, $l=0,1,2,3$.
In \eqref{18} we have
\[M_{{n}}={\frac {32}{3}}\,{\frac {\Gamma  \left( n+1/2 \right) \Gamma
 \left( n+1/8 \right) }{ \left( 16\,n-5 \right) \Gamma  \left( n+1/4
 \right) \Gamma  \left( n-1/8 \right) }} \]
and $\left\| q\right\| _{\infty } =\mathop{\max }\limits_{x\in \left[-1,1\right]} \left|q(x)\right|=1/3$. Table~3 shows that the sufficient convergence condition (\ref{34}) is fulfilled, i.e., $r_{n}<1$ for $n\ge 2$.

As in Example~\ref{ex_r1}
the computer algebra system Maple~17 was used, where the corrections to the eigenpairs \eqref{40} and \eqref{42} were computed exactly as mathematical expressions.
The approximations $\overset{m}{\lambda }_{n} $ of rank $m$ to the eigenvalues $\lambda _{n} $ for
$n=0,1,2,3,4,10$ using the FD-method of ranks $m=20$ and $m=30$ are given in Table 3.
The normalized solution of  problem (\ref{6}) with  the second boundary condition \eqref{6-1} is given by \eqref{10}.
The behavior of the corrections $\lambda _{n}^{(m)}$,  $m=0,1,...,10,20,30,$ for the eigenvalues $\lambda _{n} ,$ $n=0,1,2,3,4,$ is illustrated in Table~4.  One can observe that the method converges for $n=0,1$ too, i.e. the conditions of Theorem 3.1 are rough and can be improved. Besides we can see that the convergence of the method improves when the eigenpair number $n$ increases (see Remark 3.1).

\begin{table}
\begin{center}
\begin{tabular}{|c|c|c|c|} \hline
$n$ & $r_{n}$ & $\overset{20}{\lambda }_{n} $ & $\overset{30}{\lambda }_{n} $ \\ \hline \hline
 0 & 1.202 &  0.15067840864298071809615545712 & 0.15067840864298071808633595280 \\ \hline
 1 & 1.093 & 1.43411319565086825800583685557 & 1.43411319565086825801565633056 \\ \hline
 2 & 0.687 & 3.36697390018171063527745195043 & 3.36697390018171063527745197976 \\ \hline
 3 & 0.543 & 5.81865643939502492320620997309 & 5.81865643939502492320620997309 \\ \hline
 4 & 0.463 & 8.69578418818622939490618135946 & 8.69578418818622939490618135946 \\ \hline
10 & 0.286 & 32.6418532333094312351786092007 & 32.6418532333094312351786092007 \\ \hline
\end{tabular}
\caption{Approximations $\overset{m}{\lambda }_{n}$ of the ranks $m=20,30$ to the eigenvalues $\lambda _{n}$ and the values of the ratio $r_{n}$ of the geometric progression for $n=0,1,2,3,4,10$.}
\end{center}
\end{table}

\begin{table}
\begin{center}
\begin{tabular}{|c|c|c|c|c|c|} \hline
$m$  & $\lambda _{0}^{(m)} $ & $\lambda _{1}^{(m)}$ & $\lambda _{2}^{(m)} $ & $\lambda _{3}^{(m)} $ & $\lambda _{4}^{(m)} $ \\ \hline \hline
0 & $0.07396$ &  $1.294$ & $3.236$ &  $5.691$ & $8.569$ \\ \hline
1 & $0.08152$ & $0.1389$ & $0.1297$ &  $0.1269$ & $0.1260$ \\ \hline
2 & $-0.005109$ &  $0.001165$ & $0.001380$ &  $0.0004626$ &  $0.0002623$ \\ \hline
3 & $0.0003151$ &  $-0.0003367$ &$ 0.7655 \cdot {10}^{-5}$ &  $0.9807 \cdot {10}^{-5}$ &  $0.2132 \cdot {10}^{-5}$\\ \hline
4 & $-0.5669 \cdot {10}^{-5}$ & $0.8132 \cdot {10}^{-5}$ & $-0.2832 \cdot {10}^{-5}$  & $2.742 \cdot {10}^{-7}$ & $2.891 \cdot {10}^{-9}$ \\ \hline
5 & $-0.1551 \cdot {10}^{-5}$& $0.1676 \cdot {10}^{-5}$& $-1.164 \cdot {10}^{-7}$ & $ -9.865 \cdot {10}^{-9}$ & $ 1.484 \cdot {10}^{-9}$ \\ \hline
6 &$1.943  \cdot {10}^{-7}$ & $-1.979 \cdot {10}^{-7}$ & $4.191 \cdot {10}^{-9}$ & $-5.421 \cdot {10}^{-10}$ & $-1.394 \cdot {10}^{-11}$ \\ \hline
7 &  $-4.321 \cdot {10}^{-9}$ & $ 3.925 \cdot {10}^{-9} $& $4.053 \cdot {10}^{-10}$ & $-7.920 \cdot {10}^{-12}$ & $-1.403 \cdot {10}^{-12}$ \\ \hline
8 &$ -1.618 \cdot {10}^{-9}$ &  $1.619 \cdot {10}^{-9}$ & $-5.889 \cdot {10}^{-13}$ & $1.278 \cdot {10}^{-13}$ & $-2.140 \cdot {10}^{-14}$ \\ \hline
9 &  $2.199 \cdot {10}^{-10}$  & $-2.186 \cdot {10}^{-10}$ & $-1.253 \cdot {10}^{-12}$ & $1.466 \cdot {10}^{-14}$ & $-7.319 \cdot {10}^{-17}$ \\ \hline
10 & $-4.351 \cdot {10}^{-12}$  & $4.382 \cdot {10}^{-12}$ &$ -3.167 \cdot {10}^{-14}$ & $5.091 \cdot {10}^{-16}$ & $1.464 \cdot {10}^{-18}$ \\ \hline
20 &$ -1.116 \cdot {10}^{-20}$ & $1.116 \cdot {10}^{-20}$ & $2.576 \cdot {10}^{-26}$ & $-2.037 \cdot {10}^{-30}$ & $7.184 \cdot {10}^{-34}$ \\ \hline
30 &$ 8.072 \cdot {10}^{-30}$ & $-8.072 \cdot {10}^{-30}$ & $7.176 \cdot {10}^{-38}$ & $1.316 \cdot {10}^{-44}$ & $4.848 \cdot {10}^{-50}$ \\ \hline
\end{tabular}

\caption{Corrections $\lambda^{(m)}_{n}, n=0,1,2,3,4$ for the FD-method of the ranks $m=0,1,...,10,20,30$. }
\end{center}
\end{table}

\end{example}


\begin{example} Let us consider  problem (\ref{4}) with $\alpha =\beta =0,\, \, s =3/4$ and with the potential $q(x)=(\text{sgn}(x)+1)/2$.
The potential is not polynomial but we can use the general idea of the symbolic algorithms above representing the potential as a series  with respect to the base problem solution. Then we represent the solution $u_n^{(j+1)}(x)$ of the problems (\ref{5-1}) by a series of the same type and obtain some recurrence relations  for the series coefficients like (\ref{43})-(\ref{46}) but with $\infty$ as the upper bound for $m$. By restricting $m$ by some finite $N$ we have an approximate symbolic algorithm which we describe below in detail.

The solution of the base problem in this case is
\begin{equation} \label{48}
u_{n}^{(0)} (x)=P_{n} (x),\, \, \lambda _{n}^{(0)} =\frac{\Gamma (n+7/4)}{\Gamma (n+1/4)},
\end{equation}
where $P_{n} (x)$ is the Legendre polynomial. We look for the solution of problem (\ref{5-1}), (\ref{5-2}) in the form
\begin{equation} \label{49}
u_{n}^{(j+1)} (x)=\underset{s=0, s\ne n}{\overset{\infty}{\sum}} a_{n,s}^{(j+1)} P_{s} (x)
\end{equation}
which satisfies the  condition
\begin{equation} \label{50}
\int _{-1}^{1} u_{n}^{(j+1)} (\, x)u_{n}^{(0)} (\, x)dx=0.
\end{equation}
Let us find for the right-hand side of (\ref{5-1}) the following expansion like  (similar to \ref{49})
\begin{equation} \label{51}
F_{n}^{(j+1)} (x)=\sum _{s=0}^{\infty } f_{n,s}^{(j+1)} P_{s} (x).
\end{equation}
To find the coefficients of this expansion we use the formula
\begin{equation} \label{52}
\begin{split}
&\frac{2b_{s,t}}{2t+1}=\int _{-1}^{1} q(x)P_{s} (x)P_{t} (x)dx=\int _{0}^{1} P_{s} (x)P_{t} (x)dx \\
&{=\frac{2}{\pi (s-t)(s+t+1)} \left[A_{s,t} \sin \left(\frac{s\pi}{2}\right)\cos \left(\frac{t\pi}{2}\right)-\frac{1}{A_{s,t} } \sin \left(\frac{t\pi}{2} \right)\cos \left(\frac{s\pi}{2}\right)\right],}
\end{split}
\end{equation}
with 
$$A_{s,t} =\frac{\Gamma \left(\frac{1+s}{2} \right)\Gamma \left(1+\frac{t}{2} \right)}{\Gamma \left(\frac{1+t}{2} \right)\Gamma \left(1+\frac{s}{2} \right)} $$
(see, e.g., \cite{be}, Volume~1, Section~3.12, formula (15)). Then we obtain
\begin{equation} \label{53}
\begin{split}
&F_{n}^{(j+1)} (x)=\lambda _{n,s}^{(j+1)} P_{n} (x)+\sum _{p=1}^{j}\lambda _{n,s}^{(j-p+1)} \underset{s=0, s\ne n}{\overset{\infty}{\sum}} a_{n,s}^{(p)} P_{s} (x) \\
&{-\underset{s=0, s\ne n}{\overset{\infty}{\sum}}a_{n,s}^{(j)}  \sum _{t=0}^{\infty } b_{s,t}^{} P_{t} (x)=\left(\lambda _{n,s}^{(j+1)} -\underset{s=0, s\ne n}{\overset{\infty}{\sum}}a_{n,s}^{(j)}  b_{s,n}^{} \right)P_{n} (x)} \\
&{+\underset{s=0, s\ne n}{\overset{\infty}{\sum}} \left(\sum _{p=1}^{j}\lambda _{n,s}^{(j-p+1)}  a_{n,s}^{(p)} -\underset{t=0, t\ne n}{\overset{\infty}{\sum}}a_{n,t}^{(j)} b_{t,s}^{} \right)P_{s} (x)}.
\end{split}
\end{equation}
Now, from (\ref{5-1}), (\ref{5-2}) and (\ref{51}), we obtain the following exact symbolic algorithm
\begin{equation} \label{54}
\begin{split}
&\lambda _{n}^{(j+1)} =\underset{s=0, s\ne n}{\overset{\infty}{\sum}}a_{n,s}^{(j)}  b_{s,n}\\
&f_{n,s}^{(j+1)} =\underset{p=1}{\overset{j}{\sum}}\lambda _{n}^{(j-p+1)}  a_{n,s}^{(p)} -\underset{t=0, t\ne n}{\overset{\infty}{\sum}} a_{n,t}^{(j)} b_{t,s}, \\
&a_{n,s}^{(j+1)} =\frac{1}{\lambda _{s}^{(0)} -\lambda _{n}^{(0)} } \left(\sum _{p=1}^{j}\lambda _{n}^{(j-p+1)}  a_{n,s}^{(p)} -\underset{t=0, t\ne n}{\overset{\infty}{\sum}} a_{n,t}^{(j)} b_{t,s}\right), \\
& s=0,1,...,\infty ,\, \, \, \, s\ne n.
\end{split}
\end{equation}
The truncated sums with $N$ summands represent a practically realizable symbolic algorithm.

The computations for the case $n=0, N=64$, using Maple with 32 significant digits, resulted in
$$\overset{4}{\lambda}_0=0.73277189290359102980467600413989,$$
$$\overset{8}{\lambda}_0=0.73277298398625680196269290377257,$$
$$\overset{16}{\lambda}_0=0.73277298419141034176589978524115.$$
The correction to the eigenvalue in the 16-th step (i.e. the 16-th summand of the series) is
\[\lambda _{0}^{(16)} ={\rm 0.88694390139173379459088819979359e-13}\]
and all corrections with odd numbers are equal to zero.
It is easy to obtain analytically the corrections $\lambda _{0}^{(0)} =\frac{3\sqrt{2} }{8\pi } \Gamma ^{2} \left(\frac{3}{4} \right),\, \, u_{0}^{(0)} (x)=1,\;\lambda _{0}^{(1)} =1.$
\end{example}

{\bf{Acknowledgment.}} We want to thank the anonymous referees for their helpful remarks which contributed to improvement of the article.

\end{document}